\documentstyle[11pt]{amsart}

\title[The Lefschetz theorem for CR submanifolds]{ The Lefschetz
theorem for CR submanifolds and the nonexistence
of real analytic
Levi flat submanifolds}
\author { Lei Ni and Jon Wolfson}
\thanks{The first author was partially supported by NSF grant
DMS-0203023. The second author was partially supported by NSF
grant DMS-0104007. Both authors wish to thank the National Center for
Theoretical Sciences, National Tsing Hua University, Hsinchu, Taiwan
for the hospitality  provided during the writing of parts of this paper.}

\date{\today}

\newtheorem{thm}{Theorem}[section]

\newtheorem{cor}[thm]{Corollary}

\theoremstyle{definition}

\numberwithin{equation}{section}

\newcommand{\g}{\gamma}

\renewcommand{\l}{\lambda}

\newcommand{\n}{\nabla}

\renewcommand{\t}{\tau}
\renewcommand{\O}{\Omega}

\def\Pb{\ifmmode{\Bbb P}\else{$\Bbb P$}\fi}
\def\Z{\ifmmode{\Bbb Z}\else{$\Bbb Z$}\fi}
\def\C{\ifmmode{\Bbb C}\else{$\Bbb C$}\fi}
\def\R{\ifmmode{\Bbb R}\else{$\Bbb R$}\fi}
\def\S{\ifmmode{S^2}\else{$S^2$}\fi}

\def\S{\cal S}

\begin{document}

\maketitle

\setcounter{secnumdepth}{1}

\setcounter{section}{0}

\section{\bf Introduction}

In this  paper we study the topology of a CR submanifold with
degenerate Levi form embedded in a compact Hermitian symmetric
space. In the case that the ambient manifold is $\Pb^v$ we prove a
Lefschetz hyperplane theorem for such CR submanifolds. In addition
we relate the homotopy groups of the ambient manifold with the
homotopy groups of the CR submanifold within a range depending on
the nullity of the Levi form. The optimal results are obtained for
Levi flat submanifolds. In particular we show that a compact Levi
flat submanifold in a compact Hermitian symmetric space of complex
dimension $v$ is simply connected if the real dimension of the
Levi flat manifold is greater than $v+1$. We then apply this
result to give a proof of the nonexistence of real analytic Levi
flat submanifolds of real dimension greater than $v+1$. The proof
of the Lefschetz theorem is based on the Morse theory of paths
between a CR submanifold and  a complex submanifold. The proof of
the other restrictions on the topology of a CR submanifold is
based on the Morse theory of paths between two CR submanifolds.
The crucial ingredient in both proofs is a second variation
calculation in Schoen-Wolfson [S-W], that goes back to Frankel
[F].

To our knowledge results of this type relating the topology of a CR
submanifold with its Levi form are the first of their kind. They are
natural
generalizations of analogous results for complex
submanifolds.
On the nonexistence of Levi flat submanifolds:
Lins Neto proved the nonexistence of real analytic Levi
flat hypersurfaces in $\Pb^v$, for $v\ge 3$ [LN]; Ohsawa [O] did the
$v=2$ case,
and Siu [S] proved the nonexistence of smooth Levi
flat hypersurfaces in $\Pb^v$. Whether our nonexistence result
remains true in the smooth
category is an interesting open problem.

\bigskip

\section{\bf Definitions and the statements of the theorems}

We will use the following notion defined in [K-W]: Let $V$ be an
irreducible  compact K\"ahler manifold of complex dimension $v$
with
nonnegative holomorphic bisectional curvature. Define the symmetric
bilinear
form $H_Y(W, Z) = \langle R(Y,JY) W, JZ \rangle$, where $Y,W,Z \in
TV$. Then for any $Y \neq 0$, $H_Y$ is positive semi-definite.
Denote the null-space by ${\cal N}_Y$. Let $\ell (Y)$ be the
complementary dimension of  ${\cal N}_Y$ and define the {\it
complex positivity} $\ell$ of $V$ by
$$\ell= \inf_{Y \neq 0} \ell(Y).$$ By the classification result of
Mok [M], the above definition is independent of the choice of the
base point. An equivalent definition is used in [K-W] where $\ell$
is computed for all compact Hermitian symmetric spaces. In
particular in [K-W] it is shown:
\medskip

\noindent{${\bf \Bbb P}^v$}: \hspace{3cm} $\ell = v.$

\noindent{\bf $\mbox{Gr}(p,p + q;\C)$}: \hspace{.81cm} $\ell = p + q -1.$

\noindent{\bf $\mbox{Gr}(2,p + 2;\R)$}: \hspace{.87cm} $\ell = p -1.$

\noindent{\bf $Sp(r)/U(r)$}: \hspace{1.45cm} $\ell = r.$

\noindent{\bf $SO(2r)\!/\!U(r)$}: \hspace{1.35cm} $\ell = (r-1) + (r-2).$

\noindent{\bf $E_{6}/(Spin(10)\times T^{1})$}: \; $\ell = 11$

\noindent{\bf  $E_{7}/(E_{6}\times T^{1})$}: \hspace{1.15cm} $\ell= 17$.

\medskip

This computation goes back to [B] and [C-V]. It is also related to
the concept called {\it degree of nondegeneracy} of the
bisectional curvature in the K\"ahler geometry (cf. [S] for a
definition). In fact, if $s$ is the degree of the nondegeneracy of
$V$ one has $s+ \ell =v+1$.

Let $M$ be a compact real submanifold of real dimension $2p+1$.
For
$x \in M$, set
$T_x^{1,0}(M) = T_x^{1,0}(V) \cap (T_x(M) \otimes {\Bbb C})$. Then
$T^{1,0}(M)$
is a subbundle of  $T(M) \otimes {\Bbb C}$. We say that
$(M, T^{1,0}(M))$
is a {\it CR submanifold of $V$} with CR structure
$T^{1,0}(M)$ if the following
are satisfied:
\begin{enumerate}
\item[(i)]  ${\rm dim}_{\Bbb C} T^{1,0}(M) = p.$
\item[(ii)]
$T^{1,0}(M) \cap T^{0,1}(M) = \{0\}$, where $T^{0,1}(M)=
\overline{T^{1,0}(M)}.$
\item[(iii)] For any two sections $X_1,X_2$
of $T^{1,0}(M)$ defined on an open set $U \subset M$, the Lie bracket
$[X_1,X_2]$
is a section  of $T^{1,0}(M)$ defined on
$U$.
\end{enumerate}

Let $X_1, \dots, X_p$ be a local framing of
$T^{1,0}(M)$ on an open set $U$. Then  $\bar{X_1}, \dots,\bar{
X_p}$
is a local  framing of  $T^{0,1}(M)$ on $U$. Choose a section
$T$ of $T(M) \otimes {\Bbb C}$ so that
$X_1, \dots, X_p, \bar{X_1},
\dots,\bar{ X_p}, T$ gives a framing of $T(M) \otimes {\Bbb C}$ on
$U$. Define a hermitian matrix
$(c_{ij})$ by:
$$[X_i,\bar{X_j}] =
c_{ij}T, \;{\rm mod}\; X_1, \dots, X_p, \bar{X_1}, \dots,\bar{
X_p}.$$
The hermitian form $\t$ determined by $(c_{ij})$ is called
the {\it Levi form of the CR structure}.
The Levi matrix $(c_{ij})$
clearly depends on the choice of framing however the rank and nullity
of the associated form $\t$ are well-defined. Throughout this paper
we will call the complex
dimension of the null space of the Levi form
of $M$ the {\it nullity of $M$}.

By studying the relative geometry of a pair of submanifolds of
$\Pb^v$ consisting of a compact CR
submanifold and a compact complex
submanifold we prove the following Lefschetz theorems
for compact CR
submanifolds with degenerate Levi form.

\begin{thm}
\label{thm:Lefschetz}
Let $M$ be a compact CR submanifold of $\Pb^v$ of dimension $2p+1$.
Suppose that the nullity of $M$ is everywhere
greater than or equal to $r$, where $0 < r \leq p$. If $j \le r-1$ then
$$ \pi_j(M,M \cap \Pb^{v-1}) = 0. $$
Hence, if $j < r-1$ then the map
induced by inclusion
$$\pi_j(M \cap \Pb^{v-1}) \to \pi_j(M),$$
is an isomorphism and if $j=r-1$ this map is onto.
\end{thm}

Note that, for generic $\Pb^{v-1} \subset \Pb^{v}$, $M \cap
\Pb^{v-1}$ is a compact CR submanifold of $\Pb^{v-1}$
of dimension $2(p-1)+1$ with nullity  everywhere
greater than or equal to $r-1$, where $0 < r-1 \leq p-1$. Therefore
the result can be iterated.

We say that
$M$ is {\it Levi flat}  if the Levi form of the CR
structure vanishes. In particular, a
Levi flat submanifold admits a codimension
one foliation with complex analytic leaves.

\begin{cor}
\label{cor:lef_levi_flat}
Let $M$ be a compact Levi flat submanifold of $\Pb^v$ of
dimension $m=2p+1$.
If $j \le p-1$ then
$$ \pi_j(M,M \cap \Pb^{v-1}) = 0. $$
Hence, if $j < p-1$ then the map
induced by inclusion
$$\pi_j(M \cap \Pb^{v-1}) \to \pi_j(M),$$
is an
isomorphism and if $j=p-1$ this map is onto.
\end{cor}

For generic $\Pb^{v-1}$, $M \cap \Pb^{v-1}$ is a compact Levi flat
submanifold of $\Pb^{v-1}$
of dimension $2(p-1)+1$ and so, as above, the result can be iterated.
Combining Corollary \ref{cor:lef_levi_flat} with the
Hurewicz isomorphism theorem we get the homology version of the
Lefschetz theorem on hyperplane sections for Levi flat submanifolds
of $\Pb^v$.

\medskip

By studying the relative geometry of a pair of
compact CR submanifolds of a compact hermitian
symmetric space $V$
with dim$_{\Bbb C} V=v$  we derive  restrictions on the topology
of a compact CR submanifold with degenerate Levi form.

\begin{thm}
\label{thm:maintopology} Let $V$ be an irreducible compact K\"ahler
manifold of complex dimension $v$ with nonnegative
holomorphic bisectional curvature  of complex positivity $\ell$.
Let  $M$ be a CR submanifold
of real dimension $2p+1$ and suppose that the nullity of $M$ is
everywhere greater than or equal to
$r$, where $0 < r \leq p$. Then,
$$\pi_{j}(V,M)=0,$$
provided $j \leq \ell+1-2(v-r)$.
\end{thm}

In particular if  $M$ is Levi flat of real dimension $m$ and $ m \geq
2(v+1) - \ell$ then it follows that
$M$ is simply connected. The next
result then follows from a theorem
of Haefliger [H].

\begin{thm}
\label{thm:nonexist} Let $V$ be an irreducible compact K\"ahler
manifold of complex dimension $v$ with nonnegative
holomorphic bisectional curvature  of complex positivity $\ell$.
There are no real analytic Levi flat submanifolds of dimension
$m$ in $V$ when $m \geq 2(v+1) - \ell$.
\end{thm}

\bigskip

\section{ \bf The relative geometry of a complex submanifold and a CR
submanifold}

\setcounter{equation}{0}

Let $V$ be an irreducible compact K\"ahler manifold  of complex
dimension $v$, with nonnegative bisectional curvature and with
complex structure $J$. Let $\ell$  be the complex positivity of
$V$. Let $M$ be a real smooth compact submanifold of dimension $2p+1$
and let $N$ be a compact complex submanifold of complex dimension $n$.

We denote, by ${\O}(V;M,N) = {\O}$, the space of piecewise smooth
paths $\g:[0,1] \to V$ constrained by the requirements that
$\g(0) \in M, \g(1) \in N$. Consider the energy of a path
$$ E(\g) = \int^1_0 |\dot{\g}|^2 dt $$
as a function on ${\O}$.  It is  shown in [S-W] that
$\g$ is a critical point of $E$ if:
\begin{enumerate}
\item[(i)] $\g$ is a smooth geodesic
\item[(ii)] $\g$ is normal to $M$  at $\g(0)$ and normal to $N$ at $\g(1)$,
respectively.
\end{enumerate}

Let $W_1,W_2 \in T_{\g} {\O}$.  If $\g$ is a critical point of $E$ then
the second variation of $E$ along $\g$  is:
\begin{eqnarray*}
\frac{1}{2} E_{**}(W_1,W_2) & =  &\langle \nabla_{W_1} W_2,
\dot{\g}\rangle |^1_0
           + \int^1_0 \langle \n_{\dot{\g}} W_1, \nabla_{\dot{\g}} W_2
\rangle dt  \nonumber \\
&& - \int^1_0 \langle R(\dot{\g}, W_1) \dot{\g}, W_2 \rangle dt.
\end{eqnarray*}
where $R$ denotes the curvature tensor of $V$.

\begin{thm}
\label{thm:index1}
Let $M$ be a CR submanifold in $V$ of dimension $2p+1$. Suppose that
the nullity of $M$ is everywhere
greater than or equal to $r$, where $0 < r \leq p$.
Let $N$ be a compact complex submanifold of complex dimension $n$.
Then the index of a critical point $\g$ of $E$ is at least $\ell + 1
- (v-r) - (v-n)$,
where $\ell$ is the complex positivity of $V$.
\end{thm}

\begin{pf}
Since the null space of the Levi form of $M$ everywhere has dimension
greater than or equal to $r$, for any point $x \in M$,
locally there are ${\Bbb C}$-linearly independent vector fields
$X_1, X_2,\dots, X_{r}$ that are  sections of $T^{1,0}(M)$ and
that lie in the null space of $\t$. We can choose the $X_i$ so that
if $X_i = W_i - \sqrt{-1}JW_i$ then $W_1, \dots, W_r, JW_1, \dots,
JW_r$ are orthonormal.
Now we follow the argument of [S-W] and
estimate the index of a critical
point $\g$ of $E$.

By the second variation formula we have that
\begin{eqnarray*}
\frac{1}{2}E_{**}(W, W) &= &<\nabla _{W}W,
\dot{\g}>|^{1}_{0}+\int_0^1<\nabla_{\dot{\g}}W, \nabla_{\dot{\g}}W>\\
&\, &-\int_0^1<R(\dot{\g}, W)\dot{\g}, W>\, dt.
\end{eqnarray*}

Parallel translate $W_i$, $i=1, \dots, r$, along $\g$ to construct a
vector field
$W_i(t)$ along $\g$.  Of course, $W_i(1)$ need not be tangent to
$N$ at $\g(1)$ so $ W_i(t)$ is not necessarily an element of
$T_{\g}{\O}$. However whenever $W_i(1)\in TN$ we have:
\begin{equation}
\frac{1}{2}E_{**}(W_i, W_i)=<\nabla _{W_i}W_i,
\dot{\g}>|^{1}_{0}-\int_0^1<R(\dot{\g}, W_i)\dot{\g}, W_i>\, dt.
\end{equation}
Similarly, parallel translate  $JW_{i}$ along $\g$ to get:
\begin{equation}
\frac{1}{2}E_{**}(JW_i, JW_i)=<\nabla _{JW_i}JW_i,
\dot{\g}>|^{1}_{0}-\int_0^1<R(\dot{\g}, JW_i)\dot{\g}, JW_i>\, dt.
\end{equation}
Observe that,
\begin{eqnarray*}
<\nabla_{JW_i} JW_i, \dot{\g}>& =& <J(\nabla_{JW_i}W_i), \dot{\g}>\\
& = & -<\nabla_{JW_i}W_i, J \dot{\g}>\\
& = & -<\nabla_{W_i}JW_i, J\dot{\g}>+<[JW_i, W_i], J\dot{\g}>\\
& = & -<\nabla_{W_i}W_i, \dot{\g}>+<[JW_i, W_i], J\dot{\g}>,
\end{eqnarray*}
At $\g(0)$ the vector $\dot{\g}$ is normal to $M$. Thus at  $\g(0)$
the vector $J\dot{\g}$ is either  normal to $M$ or, if tangent to $M$,
then it is normal to $T^{1,0}(M) \oplus T^{0,1}(M)$. By assumption
$X_i(0)$  lies in the null space of the Levi form, so at $\g(0)$,
$<[JW_i, W_i],
J\dot{\g}>=0$.
To show that $<[JW_i, W_i], J\dot{\g}>=0$ at  $\g(1)$
we use a simple argument of Frankel [F]. Let $R$ be an
analytic curve that passes through $\g(1)$ with tangent vectors $W_i,
JW_i$ at $\g(1)$ and
extend $W_i, JW_i$ to be tangent vector fields to $R$ in a neighborhood
$U$ of $\g(1)$.
Then on $U$, $[JW_i,W_i]$ and $J[JW_i,W_i]$ are tangent to $R$ and
so, at $\g(1)$, they
are  perpendicular to $\dot{\g}(1)$.

Therefore adding (2.1) and (2.2) we have
\begin{eqnarray*}
\frac{1}{2}E_{**}(W_i, W_i) &+ &\frac{1}{2}E_{**}(JW_i, JW_i)=\\
&-&\int_0^1\left(<R(\dot{\g}, W_i)\dot{\g}, W_i>+<R(\dot{\g},
JW_i)\dot{\g}, JW_i>\right)\, dt.
\end{eqnarray*}
Using the symmetries of the curvature tensor we have:
\begin{equation}
\langle R(\dot{\g},W_i)\dot{\g},W_i \rangle + \langle
R(\dot{\g},JW_i)\dot{\g},JW_i \rangle = \langle R(\dot{\g},J
\dot{\g})W_i ,JW_i \rangle.
\label{equ:curv}
\end{equation}
This expression is the holomorphic bisectional curvature of the complex
lines $\dot{\g} \wedge J \dot{\g}$ and $W_i \wedge JW_i$.

Note that the vectors $W_i(1), JW_i(1)$ are perpendicular to both
$\dot{\g}(1)$ and $J\dot{\g}(1)$.  Thus the vector space
$$ S = \mbox{span}\{ W_1(1), \dots, W_{r}(1), JW_1(1), \dots, JW_{r}(1)\}$$
is a complex $r$-dimensional space lying in a complex $(v-1)$-dimensional
subspace of $T_{\g(1)}V$.  It follows that the subspace $S \cap
T_{\g(1)}N$
has complex dimension at least equal to
$$ r + n - (v-1). $$
Moreover, the vector fields $\{W,JW\}$ with $W(1), JW(1) \in S \cap T_{\g(1)}N$
are parallel and lie in $T_{\g} {\O}$. Therefore
under the assumption that the positivity of $V$ is $\ell$ it
follows that every nontrivial critical point $\g$ of $E$ on
$\Omega$ has index $\lambda \ge \lambda_0$ with,
$$\l_0 =  r + n - (v-1) - (v -\ell).$$
\end{pf}

Using the relation between Morse theory on the path space $\Omega$
and the topology of $(V; M,N)$ (see [S-W] or [K-W] for details) we have:

\begin{thm}
\label{thm:top-pair}
Let $V$ be a compact K\"ahler manifold of non-negative holomorphic
bisectional curvature, of complex dimension $v$ and with complex
positivity $\ell$.
Let $M$ be a compact CR submanifold of $V$ of
dimension $m=2k+1$ and suppose that the nullity of $M$ is everywhere
greater than or equal to $r$, where $0 < r \leq k$. Let $N$ be a
compact complex submanifold of complex dimension $n$. Then
the homomorphism induced by the inclusion
$$ \imath_*: \pi_j(M,M \cap N) \to \pi_j(V,N) $$
is an isomorphism for $j \le n+r-v-(v-\ell)$ and is a surjection for
$j = n+r-v-(v-\ell) + 1$ .
\end{thm}

\bigskip

The following result is proved in [S-W] and [K-W] using techniques
analogous to those used above.

\begin{thm}
\label{thm:SW}
Let $V$ be a compact K\"ahler manifold of non-negative holomorphic
bisectional curvature, of complex dimension $v$ and with complex
positivity $\ell$.
Let $N$ be a compact complex submanifold of
complex dimension $n$.
If
$$j \le 2n-v-(v-\ell)+1$$
then
$$ \pi_j(V,N) = 0. $$
\end{thm}

\bigskip

\begin{cor}
\label{cor:general_lef}
Under the same hypothesis as in Theorem \ref{thm:top-pair}, if
$$j \le \min(2n-v-(v-\ell)+1, ~ n+r -v-(v-\ell))$$
then
$$ \pi_j(M,M \cap N) = 0. $$
\end{cor}

\begin{pf}
Follows from Theorem \ref{thm:SW} and  Theorem \ref{thm:top-pair}.
\end{pf}

\bigskip

Applying Corollay \ref{cor:general_lef} to the special case when  $V
= \Pb^v$ with
the Fubini-Study metric and $N= \Pb^{v-1}$ is a hyperplane we have
the Lefschetz hyperplane
theorem  for compact CR submanifolds with
degenerate Levi form (Theorem \ref{thm:Lefschetz})
and for Levi flat
submanifolds (Corollary \ref{cor:lef_levi_flat}).

\bigskip

\section{ \bf The relative geometry of a pair of CR
submanifolds}

Let $V$ be an irreducible compact K\"ahler manifold
of complex
dimension $v$, with nonnegative bisectional curvature and with
complex structure $J$. Let $\ell$  be the complex positivity of
$V$. Let $M$ and $N$ be a real smooth compact submanifolds of
dimensions $m=2p+1$
and $n=2q+1$, respectively.

We denote, by ${\O}(V;M,N) = {\O}$, the space of piecewise smooth
paths $\g:[0,1] \to V$ constrained by the requirements that
$\g(0) \in M, \g(1) \in N$. Consider the energy of a path
$$ E(\g) = \int^1_0 |\dot{\g}|^2 dt $$
as a function on ${\O}$.  It is  shown in [S-W] that
$\g$ is a critical point of $E$ if:
\begin{enumerate}
\item[(i)] $\g$ is a smooth geodesic
\item[(ii)] $\g$ is normal to $M$  at $\g(0)$ and normal to $N$ at $\g(1)$,
respectively.
\end{enumerate}

Let $W_1,W_2 \in T_{\g} {\O}$.  If $\g$ is a critical point of $E$ then
the second variation of $E$ along $\g$  is:
\begin{eqnarray*}
\frac{1}{2} E_{**}(W_1,W_2) & =  &\langle \nabla_{W_1} W_2,
\dot{\g}\rangle |^1_0
           + \int^1_0 \langle \n_{\dot{\g}} W_1, \nabla_{\dot{\g}} W_2
\rangle dt  \nonumber \\
&& - \int^1_0 \langle R(\dot{\g}, W_1) \dot{\g}, W_2 \rangle dt.
\end{eqnarray*}
where $R$ denotes the curvature tensor of $V$.

\begin{thm}
\label{thm:index1}
Let $M$ and $N$ be  CR submanifolds of $V$ of dimensions $2p+1$ and
$2q+1$, respectively.
Suppose that the nullity of $M$ and $N$  are
everywhere
greater than or equal to $r$ and $s$, respectively, where $0 < r \leq
p$ and $0 < s \leq q$.
Then the index  of a critical point $\g$ of $E$ is at least $\ell + 1 - (v-r) -
(v-s)$,
where $\ell$ is the complex positivity of $V$.
\end{thm}

\begin{pf}
Since the null space of the Levi form of $M$ everywhere has dimension
greater than or equal to $r$, for any point $x \in M$,
locally there are ${\Bbb C}$-linearly independent vector fields
$X_1, X_2,\dots, X_{r}$ that are  sections of $T^{1,0}(M)$ and that
lie in the null space of the Levi form of $M$. We can choose the $X_i$
so that if $X_i = W_i - \sqrt{-1}JW_i$
then $W_1, \dots, W_r, JW_1, \dots,
JW_r$ are orthonormal.
Now we follow the argument of [S-W] and
estimate the index of a critical point $\g$ of $E$.

By the second variation formula we have that
\begin{eqnarray*}
\frac{1}{2}E_{**}(W, W) &= &<\nabla _{W}W,
\dot{\g}>|^{1}_{0}+\int_0^1<\nabla_{\dot{\g}}W, \nabla_{\dot{\g}}W>\\
&\, &-\int_0^1<R(\dot{\g}, W)\dot{\g}, W>\, dt.
\end{eqnarray*}

Parallel translate $W_i$, $i=1, \dots, r$, along $\g$ to construct a
vector field
$W_i(t)$ along $\g$.  Of course, $W_i(1)$ need not be tangent to
$N$ at $\g(1)$ so $ W_i(t)$ is not necessarily an element of
$T_{\g}{\O}$. However whenever $W_i(1)\in TN$ we have:
\begin{equation}
\frac{1}{2}E_{**}(W_i, W_i)=<\nabla _{W_i}W_i,
\dot{\g}>|^{1}_{0}-\int_0^1<R(\dot{\g}, W_i)\dot{\g}, W_i>\, dt.
\end{equation}
Similarly, parallel translate  $JW_{i}$ along $\g$ to
get:
\begin{equation}
\frac{1}{2}E_{**}(JW_i, JW_i)=<\nabla _{JW_i}JW_i,
\dot{\g}>|^{1}_{0}-\int_0^1<R(\dot{\g}, JW_i)\dot{\g}, JW_i>\, dt.
\end{equation}
Now assume that $X_i(1)= W_i(1) - \sqrt{-1}JW_i(1)$
lies in the null space of the Levi form of $N$ at $\g(1)$.
Observe that,
\begin{eqnarray*}
<\nabla_{JW_i} JW_i, \dot{\g}>& =& <J(\nabla_{JW_i}W_i), \dot{\g}>\\
& = & -<\nabla_{JW_i}W_i, J \dot{\g}>\\
& = & -<\nabla_{W_i}JW_i, J\dot{\g}>+<[JW_i, W_i], J\dot{\g}>\\
& = & -<\nabla_{W_i}W_i, \dot{\g}>+<[JW_i, W_i], J\dot{\g}>,
\end{eqnarray*}
At $\g(0)$ the vector $\dot{\g}$ is normal to $M$. Thus at  $\g(0)$
the vector $J\dot{\g}$ is either  normal to $M$ or, if tangent to $M$,
then it is normal to $T^{1,0}(M) \oplus T^{0,1}(M)$. By assumption
$X_i(0)$  lies in the null space of the Levi form, so at $\g(0)$,
$<[JW_i, W_i],
J\dot{\g}>=0$.
To show that $<[JW_i, W_i], J\dot{\g}>=0$ at $\g(1)$
we use the assumption that $X_i(1)$ lies in the null space of the
Levi form
of $N$ at $\g(1)$ and apply the previous reasoning.
Therefore adding (2.1) and (2.2) we have
\begin{eqnarray*}
\frac{1}{2}E_{**}(W_i, W_i) &+ &\frac{1}{2}E_{**}(JW_i, JW_i)=\\
&-&\int_0^1\left(<R(\dot{\g}, W_i)\dot{\g}, W_i>+<R(\dot{\g},
JW_i)\dot{\g}, JW_i>\right)\, dt.
\end{eqnarray*}
Using the symmetries of the curvature tensor we have:
\begin{equation}
\langle R(\dot{\g},W_i)\dot{\g},W_i \rangle + \langle
R(\dot{\g},JW_i)\dot{\g},JW_i \rangle = \langle R(\dot{\g},J
\dot{\g})W_i ,JW_i \rangle.
\label{equ:curv}
\end{equation}
This expression is the holomorphic bisectional curvature of the complex
lines $\dot{\g} \wedge J \dot{\g}$ and $W_i \wedge JW_i$.

Note that the vectors $W_i(1), JW_i(1)$ are perpendicular to both
$\dot{\g}(1)$ and $J\dot{\g}(1)$.  Thus the vector space
$$ S = \mbox{span}\{ W_1(1), \dots, W_{r}(1), JW_1(1), \dots, JW_{r}(1)\}$$
is a complex $r$-dimensional space lying in a complex $(v-1)$-dimensional
subspace of $T_{\g(1)}V$.  Denote the null space of the Levi form of
$N$ at $x \in N$
by ${\cal Z}_x(TN)$. Then the subspace $S \cap {\cal
Z}_{\g(1)}(TN)$
has complex dimension at least equal to
$$ r + s - (v-1). $$
Moreover, the vector fields $\{W,JW\}$ with $W(1), JW(1) \in S \cap
{\cal Z}_{\g(1)}(TN)$
are parallel and lie in $T_{\g} {\O}$. Therefore
under the assumption that the positivity of $V$ is $\ell$ it
follows that every nontrivial critical point $\g$ of $E$ on
$\Omega$ has index $\lambda \ge \lambda_0$ with,
$$\l_0 =  r + s - (v-1) - (v -\ell).$$
\end{pf}

Using the relation between Morse theory on the path space $\Omega$
and the topology of $(V; M,N)$ (see [S-W] or [K-W] for details) we have:

\begin{thm}
\label{thm:top-pair2}
Let $V$ be a compact K\"ahler manifold of non-negative holomorphic
bisectional curvature, of complex dimension $v$ and with complex
positivity $\ell$.
Let $M,N$ be a compact CR submanifolds of $V$ of
dimensions $2p+1$ and $2q+1$, respectively and
suppose that the nullity of $M$ and $N$ are everywhere
greater than or equal to $r$ and $s$, respectively, where $0 < r \leq
p$ and $0 < s \leq q$.
Then the homomorphism induced by the inclusion
$$ \imath_*: \pi_j(M,M \cap N) \to \pi_j(V,N) $$
is an isomorphism for $j \le r+s-v-(v-\ell)$ and is a surjection for
$j = r+s-v-(v-\ell) + 1$ .
\end{thm}

\bigskip

Applying this result to the case $M=N$ we have:

\begin{cor}
Let $M$ be a CR submanifold of $V$ as in Theorem \ref{thm:top-pair2} .
Then
$$\pi_{j}(V,M)=0,$$
provided $j \leq \ell+1-2(v-r)$, where $\ell$ is the complex positivity of $V$.
\end{cor}

Using the exact homotopy sequence of the pair $(V,M)$ and that $\pi_1(V)=0$ we
conclude:

\begin{cor}
\label{cor:CRfundgroup} Let $V$ and $M$ be as above. If
$2 r \geq 2v+1-\ell$  then $\pi_1(M)=\{0\}$.
\end{cor}

\begin{cor}
\label{cor:fundgroup} Let $V$  be as above and suppose that $M$ is
Levi flat. If {\rm dim} $M \geq 2(v+1)-\ell$ then $\pi_1(M)=\{0\}$.
\end{cor}

Recall the well known result of A. Haefliger [H].
\begin{thm}
\label{thm:haef}
There is no real analytic, codimension one foliation of a compact
simply-connected manifold.
\end{thm}

Corollary \ref{cor:fundgroup} and Theorem \ref{thm:haef} imply
Theorem \ref{thm:nonexist}.

\bigskip

Combining the topological restrictions imposed on Levi flat
submanifolds by Theorem \ref{thm:top-pair2}
and its corollaries with codimension one foliation theory leads to
many partial nonexistence results
such as:

\begin{thm}
Let $V$ be an irreducible compact K\"ahler
manifold of complex dimension $v$ with complex positivity $\ell$.
If  $m \geq 2(v+1) - \ell$ then there are no smooth Levi flat
submanifolds of dimension
$m$ in $V$ with a compact leaf.
\end{thm}

We leave the formulation of similar results to the reader.

\end{document}